\documentstyle[11pt,mst-stylefile,harvard,Bbb11a]{article}
\catcode`@=12
\newcommand{\bbZ}{{\Bbb Z}}
\newcommand{\bbR}{{\Bbb R}}

\newcommand{\bbL}{{\Bbb L}}

\renewcommand{\cite}{\citeyear}

\begin{document}

\title{Semi-additive functionals and cocycles \\  in the context of self-similarity
\thanks{ This research was partially supported by the NSF grant DMS-0102410 at Boston University. The second author
would also like to thank SAMSI and the University of North Carolina at Chapel Hill for
their hospitality.}
 }

\author{
 Vladas Pipiras
\\  University of North Carolina at Chapel Hill
\and
 Murad S.Taqqu
\\  Boston University
}

\bibliographystyle{agsm}

\maketitle


\begin{abstract}
 \noindent

Self-similar symmetric $\alpha$-stable, $\alpha\in(0,2)$, mixed moving averages can be
related to nonsingular flows. By using this relation and the structure of the underlying
flows, one can decompose self-similar mixed moving averages into distinct classes and
then examine the processes in each of these classes separately. The relation between
processes and flows involves semi-additive functionals. We establish a general result
about semi-additive functionals related to cocycles, and identify the presence of a new
semi-additive functional in the relation between processes and flows. This new functional
is useful for finding the kernel function of self-similar mixed moving averages generated
by a given flow. It also sheds new light on previous results on the subject.
\end{abstract}

\section{Introduction}
\label{s:preliminar}

\noindent A symmetric $\alpha$-stable ($S\alpha S$, in short), $\alpha\in (0,2)$, process
$\{X(t)\}_{t\in\bbR}$ is called a (stationary increments) {\it mixed moving average} if
it has an integral representation of the form
\begin{equation}\label{e:mma}
\{X_\alpha(t)\}_{t\in\bbR} \stackrel{d}{=}\left\{ \int_X\int_\bbR \Big(G(x,t+u) -
G(x,u)\Big) M_\alpha(dx,du)\right\}_{t\in\bbR}
\end{equation}
or, equivalently, if the characteristic function of the process $X_{\alpha}$ is given by
\begin{equation}\label{e:char-func-mma}
E\exp\Big\{i\sum_{k=1}^n \theta_k X_\alpha(t_k) \Big\} = \exp\biggl\{-\int_X \int_\bbR
\Big| \sum_{k=1}^n \theta_k G_{t_k}(x,u) \Big|^\alpha \mu(dx)du \biggr\}.
\end{equation}
In these relations, $\stackrel{d}{=}$ denotes the equality in the sense of the
finite-dimensional distributions, $(X,{\cal X},\mu)$ is a standard Lebesgue space,
$M_\alpha$ is a $S\alpha S$ random measure on $X\times \bbR$ with control measure
$\mu(dx)du$ and $G:X\times\bbR\mapsto \bbR$ is a deterministic function, called a {\it
kernel (function)}, such that the ``time increment''
\begin{equation}\label{e:G_t}
    G_t(x,u) = G(x,t+u) - G(x,u),\quad x\in X,u\in\bbR,
\end{equation}
satisfies $\{G_t\}_{t\in\bbR}\subset L^\alpha(X\times \bbR,\mu(dx)du)$. For more
information on $S\alpha S$ random measures, processes and integral representations of the
type \refeq{mma}, see for example Samorodnitsky and Taqqu
\cite{samorodnitsky:taqqu:1994book}. Observe that, in view of (\ref{e:char-func-mma}),
mixed moving averages have always stationary increments. They are symmetric
$\alpha$-stable stationary increments processes related to dissipative flows in the sense
of Surgailis, Rosi{\'n}ski, Mandrekar and Cambanis
\cite{surgailis:rosinski:mandrekar:cambanis:1998}, and form an important subclass of all
symmetric $\alpha$-stable processes with stationary increments. They also can be viewed
as stationary increments extensions of {\it stationary} mixed moving averages of Rosi{\'
n}ski \cite{rosinski:1995}.

To avoid trivialities, we will assume that
\begin{equation}\label{e:full-support}
\mbox{supp}\left\{ G_t(x,u),t\in\bbR\right\} = X\times\bbR \quad \mbox{a.e.}\ \mu(dx)du
\end{equation}
holds. By support $\mbox{supp} \{G_t, t \in \bbR\}$ we mean a minimal (a.e.) set $A
\subset X\times \bbR$ such that $m \{G_t(x,u) \neq 0, (x,u) \notin A\}=0$ for every $t
\in \bbR$, where $dm = \mu(dx)du$.

We shall focus on $S\alpha S$, $\alpha\in (0,2)$, mixed moving averages
$\{X_\alpha(t)\}_{t \in \bbR}$ which are, in addition, {\it self-similar} with a
self-similarity parameter $H>0$, that is, for any $c>0,$
\begin{equation} \label{e:ss}
\{ X_{\alpha}(ct)\}_{t \in \bbR} \stackrel{d}{=} \{ c^H
X_{\alpha}(t)\}_{t \in \bbR}.
\end{equation}
Such processes are of interest because their stationary increments can be used as models
for strongly dependent $S\alpha S$ time series. Observe, however, that self-similarity
imposes further restrictions on the kind of kernel functions $G$ which can be used in the
representation (\ref{e:mma}). Self-similar $S\alpha S$ mixed moving averages are studied
in Pipiras and Taqqu
\cite{pipiras:taqqu:2002d,pipiras:taqqu:2002s,pipiras:taqqu:2003di,pipiras:taqqu:2003id,pipiras:taqqu:2003in}.
One can relate these processes to nonsingular flows. By using this relation and the
structure of these flows, one can decompose self-similar mixed moving averages into
distinct classes and then examine processes in each of the classes separately. The
relation between self-similar mixed moving averages and flows established by Pipiras and
Taqqu \cite{pipiras:taqqu:2002d} (see Definition 5.1 and Proposition 5.1) is as follows.

\begin{definition}\label{d:ss-mma} (Pipiras and Taqqu \cite{pipiras:taqqu:2002d})
A $S\alpha S$, $\alpha\in(0,2)$, self-similar process $X_\alpha$ having a mixed moving
average representation (\ref{e:mma}) is said to be {\it generated by a nonsingular
measurable flow} $\{\psi_c\}_{c>0}$ on $(X,{\cal X},\mu)$ if, for all $c>0$,
\begin{equation}\label{e:generated}
c^{-(H-1/\alpha)} G(x,cu) =  b_c(x) \left\{{d(\mu\circ \psi_c) \over d\mu}
(x)\right\}^{1/\alpha} G \Big(\psi_c(x), u + g_c(x) \Big) + j_c(x),\quad \mbox{a.e.}\
\mu(dx)du,
\end{equation}
where $\{b_c\}_{c>0}$ is a cocycle for the flow $\{\psi_c\}_{c>0}$ taking values in
$\{-1,1\}$, $\{g_c\}_{c>0}$ is a semi-additive functional for the flow $\{\psi_c\}_{c>0}$
and $j_c(x):(0, \infty) \times X \to \bbR$ is some function, and also if
(\ref{e:full-support}) holds.
\end{definition}

We shall now define the various terms used in Definition \ref{d:ss-mma}. A
(multiplicative) {\it flow} $\{\psi_c\}_{c>0}$ is a collection of deterministic maps
$\psi_c: X \to X$ such that
\begin{equation}\label{e:flow}
\psi_{c_1c_2}(x) = \psi_{c_1}(\psi_{c_2}(x)),\ {\mbox {for all}} \ c_1,c_2>0,x\in X,
\end{equation}
and $\psi_1(x)=x$, for all $x\in X$. A flow is {\it nonsingular} if $\mu(A)=0$ implies
$\mu(\psi_c^{-1}(A))=0$ for any $c>0$ and $A \in {\cal X}$. It is {\it measurable} if the
map $\psi_c(x):(0,\infty)\times X\mapsto X$ is measurable. A { \it cocycle}
$\{b_c\}_{c>0}$ for the flow $\{\psi_c\}_{c>0}$ is a measurable map
$b_c(x):(0,\infty)\times X\mapsto Y$ satisfying the relation
\begin{equation}\label{e:cocycle}
b_{c_1c_2}(x) = b_{c_1}(x)b_{c_2}(\psi_{c_1}(x)),\ {\mbox {for all}} \ c_1,c_2>0,x\in X.
\end{equation}
In our context, we will have either $Y=\{-1,1\}$, $Y=(0,\infty)$ or $Y = \bbR\setminus
\{0\}$. A {\it semi-additive functional} $\{g_c\}_{c>0}$ for the flow $\{\psi_c\}_{c>0}$
is a measurable function $g_c(x):(0,\infty)\times X\mapsto \bbR$ satisfying
\begin{equation}\label{e:I-semi-add}
g_{c_1c_2}(x) = \frac{g_{c_1}(x)}{c_2} + g_{c_2}(\psi_{c_1}(x)),\ {\mbox {for all}} \
c_1,c_2>0,x\in X.
\end{equation}

All quantities entering in (\ref{e:generated}) thus obey a specific relation with the
exception of the function $j_c(x)$. One of the goals of this work is to show that,
without loss of generality, this function also obeys a special relation, namely,
\begin{equation}\label{e:II-semi-add}
j_{c_1c_2}(x) = c_2^{-(H-1/\alpha)} j_{c_1}(x) + b_{c_1}(x) \left\{\frac{d(\mu\circ
\psi_{c_1})}{d\mu}(x) \right\}^{1/\alpha} j_{c_2}(\phi_{c_1}(x)),\ \mbox{for all} \
c_1,c_2>0,x\in X.
\end{equation}
It is quite easy to show that $\{j_c\}_{c>0}$ satisfies (\ref{e:II-semi-add}) a.e.\
$\mu(dx)$, for all $c_1,c_2>0$ (see the proof of Theorem \ref{t:2-semi-additive}), in
which case $\{j_c\}_{c>0}$ is called {\it an almost {\it 2}-semi-additive functional}.
Having (\ref{e:II-semi-add}) for a.e.\ $\mu(dx)$, however, is often not enough or
difficult to solve for $\{j_c\}_{c>0}$ given a specific flow $\{\psi_c\}_{c>0}$.
Therefore, we need to show that $\{j_c\}_{c>0}$ has a version (an a.e.\ modification for
fixed $c$) satisfying (\ref{e:II-semi-add}) for all $x\in X$, which entails the
specification of an adequate version of the Radon-Nikodym derivative $d(\mu\circ
\psi_{c})/d\mu$. A functional $\{j_c\}_{c>0}$ satisfying (\ref{e:II-semi-add}) will be
called {\it 2-semi-additive functional} (Definition \ref{d:1-2-semi-additive}).

Since $g_c$ plays a role parallel to $j_c$, a functional $\{g_c\}_{c>0}$ satisfying
(\ref{e:I-semi-add}) will be called {\it 1-semi-additive functional}. {\it
2}-semi-additive and other functionals are important when solving the equation
(\ref{e:generated}) for the kernel function $G$ related to some specific flows (see, for
example, Pipiras and Taqqu \cite{pipiras:taqqu:2003in}, Theorem 3.1 and its proof). {\it
2}-semi-additive functionals also shed light on the structure of our previous results
(see the discussion following the proof of Theorem \ref{t:2-semi-additive} below and the
remark following Example \ref{ex:1-semi-additive-ident-dissip} below).

The proof of the existence of a version which makes (\ref{e:II-semi-add}) valid for all
$x\in X$, is involved. We establish this result in a more general context, namely, for
the so-called {\it semi-additive functionals related to a cocycle}. After a simple
transformation, {\it 1}- and {\it 2}-semi-additive functionals are examples of
semi-additive functionals related to particular cocycles (see Examples
\ref{ex:semi-additive-cocycle2} and \ref{ex:semi-additive-cocycle1}).

 We show not only
the existence of a version but also obtain an expression for  semi-additive
functionals related to cocycles. We then use this expression to characterize {\it 1}- and
{\it 2}-semi-additive functionals associated with cyclic flows.

The paper is organized as follows. Section \ref{s:functionals-cocycle} contains results
on the semi-additive functionals related to cocycles. In Section \ref{s:new-functional},
we prove that the function $j_c(x)$ in (\ref{e:generated}) can be taken as a {\it
2}-semi-additive functional satisfying (\ref{e:II-semi-add}). We provide a number of
examples in Section \ref{s:examples} which illustrate how one can derive the
semi-additive functionals. Semi-additive functionals associated with cyclic flows are
studied in Section \ref{s:functionals-cyclic}.


\section{Semi-additive functionals related to cocycles}
\label{s:functionals-cocycle}

For notational convenience, we will work here with additive flows and related
functionals. Let $\{\phi_t\}_{t\in\bbR}$ be a measurable (additive) flow on a standard
Lebesgue space $(X,{\cal X},\mu)$ satisfying $$\phi_{t_1+t_2}(x) =
\phi_{t_1}(\phi_{t_2}(x))$$ for all $t_1,t_2\in \bbR$, $x\in X$, and $\phi_0(x) = x$ for
all $x\in X$. Let $\{A_t\}_{t\in\bbR}$ be a cocycle for the flow $\{\phi_t\}_{t\in\bbR}$
taking values in $\bbR\setminus \{0\}$, that is, a functional satisfying
\begin{equation}\label{e:coc-add}
A_{t_1+t_2}(x) = A_{t_1}(x)A_{t_2}(\phi_{t_1}(x))
\end{equation}
for all $t_1,t_2\in \bbR$, $x\in X$. Observe that
\begin{equation}\label{e:a-zero}
    A_0(x) = 1
\end{equation}
since (\ref{e:coc-add}) implies $A_0(x) = A_0(x)A_0(\phi_0(x)) = (A_0(x))^2$ by using
$\phi_0(x)=x$.

\begin{definition}\label{d:semi-additive-cocycle}
A measurable map $F_t(x):\bbR \times X\mapsto \bbR$ is called a {\it semi-additive
functional related to a cocycle} $\{A_t\}_{t\in\bbR}$ (and a flow
$\{\phi_t\}_{t\in\bbR}$) if, for all $t_1,t_2\in\bbR$, $x\in X$,
\begin{equation}\label{e:semi-additive-cocycle}
F_{t_1+t_2}(x) = F_{t_1}(x) + A_{t_1}(x) F_{t_2}(\phi_{t_1}(x)).
\end{equation}
When $F_t(x)$ satisfies (\ref{e:semi-additive-cocycle}) a.e.\ $\mu(dx)$, for all
$t_1,t_2\in\bbR$, it is called an {\it almost semi-additive functional related to a
cocycle}.
\end{definition}

\noindent{\bf Remark.} In the context of multiplicative flows $\{\psi_c\}_{c>0}$, a
semi-additive functional $\{J_c\}_{c>0}$ related to a cocycle $\{B_c\}_{c>0}$ is defined
by the relation
\begin{equation}\label{e:semi-additive-cocycle-m}
J_{c_1c_2}(x) = J_{c_1}(x) + B_{c_1}(x) J_{c_2}(\psi_{c_1}(x)),
\end{equation}
for all $c_1,c_2>0$, $x\in X$.

\medskip

\noindent{\bf Notation.} We write $\{\phi_t\}_{t\in\bbR}$, $\{A_t\}_{t\in\bbR}$ and
$\{F_t\}_{t\in\bbR}$ in the context of additive flows and $\{\psi_c\}_{c>0}$,
$\{B_c\}_{c>0}$ and $\{J_c\}_{c>0}$ in the context of multiplicative flows to denote
respectively flows, cocycles and semi-additive functionals.

\medskip

In the next example, we introduce a large class of semi-additive functionals related to
cocycles.

\begin{example}\label{ex:semi-additive-cocycle-common}
Let $\{A_t\}_{t\in\bbR}$ be a cocycle for the flow $\{\phi_t\}_{t\in\bbR}$ taking values
in $\bbR\setminus \{0\}$. Let also $F:X\mapsto \bbR$ be a function and set
\begin{equation}\label{e:semi-additive-cocycle-common}
F_t(x) = A_t(x)F(\phi_t(x)) - F(x).
\end{equation}
By using the cocycle equation (\ref{e:coc-add}), we have for $t_1,t_2\in\bbR$,
\begin{eqnarray*}
  F_{t_1+t_2}(x) &=& A_{t_1}(x)A_{t_2}(\phi_{t_1}(x)) F(\phi_{t_2}(\phi_{t_1}(x))) - F(x) \\
   &=& A_{t_1}(x) \Big( A_{t_2}(\phi_{t_1}(x)) F(\phi_{t_2}(\phi_{t_1}(x))) -
F(\phi_{t_1}(x))\Big)  \\
& & \quad \quad \quad +\ A_{t_1}(x)F(\phi_{t_1}(x)) - F(x) \\
   &=& A_{t_1}(x) F_{t_2}(\phi_{t_1}(x)) + F_{t_1}(x).
\end{eqnarray*}
Therefore, in view of (\ref{e:semi-additive-cocycle}), the functional
$\{F_t\}_{t\in\bbR}$ in (\ref{e:semi-additive-cocycle-common}) is a semi-additive
functional related to the cocycle $\{A_t\}_{t\in\bbR}$.
\end{example}

Viewed from a different angle, (\ref{e:semi-additive-cocycle-common}) is a particular
solution to the equation (\ref{e:semi-additive-cocycle}). It is not necessarily the
unique solution. We will encounter in Section \ref{s:examples} a number of particular
equations of the form (\ref{e:semi-additive-cocycle}), for which we derive the general
solutions.

We now show that an almost semi-additive functional related to a cocycle has a version
which is a semi-additive functional related to a cocycle. We say that $\{f_t\}_{t\in
T}\subset L^0(S,{\cal S},m)$ is a {\it version} of $\{\widetilde f_t\}_{t\in T}\subset
L^0(S,{\cal S},m)$, where $T$ is an arbitrary index set, if $m(f_t\neq \widetilde f_t) =
0$ for all $t\in T$. The proof of this result uses the notion of a special flow
$\widetilde \phi_t(y,u)$ . Informally, the flow $\widetilde \phi_t(y,u)$ is defined on
the set of points
$$\Omega = \{(y,u):0\le u<r(y),y\in Y\} =
Y\times [0,r(\cdot)),$$ where $r(y)$ is a positive function. Plotting $(y,u)$ in two
dimensions, we can view the flow $\widetilde \phi_t$ as moving up vertically at constant
speed till it reaches the level $r(y)$, and then jumps back to a point $(y',0)$ before it
renews its vertical climb, this time from the point $y'$.  Thus, if one focuses only on
the horizontal $Y$ axis, the flow starting at $y$ moves to $y' = Vy$, then to
$V^2y,\dots,V^ny,\dots$ Since the flow $\widetilde \phi_t$ moves constantly, observe that
it has no fixed points.

We shall apply Theorem 3.1 of Kubo \cite{kubo:1969} to define $\widetilde \phi_t$
formally. According to that theorem, any (measurable, nonsingular) flow
$\{\phi_t\}_{t\in\bbR}$ without fixed points on a standard Lebesgue space is
null-isomorphic (mod 0) to some {\it special flow} $\{\widetilde \phi_t\}_{t\in\bbR}$
defined on the space $(Y\times [0,r(\cdot)),{\cal Y}\otimes {\cal B}([0,r(\cdot))),
\tau(dy)du)$ by
\begin{equation}\label{e:special-representation}
\widetilde \phi_t(y,u) = (V^n y,t+u-r_n(y)), \quad \mbox{for}\ 0\leq u+t-r_n(y) <
r(V^ny),
\end{equation}
where $(Y,{\cal Y},\tau)$ is a standard Lebesgue space, $V$ is a null-isomorphism of $Y$
onto itself, $r$ is a positive measurable function on $Y$ satisfying
$\sum_{k=-\infty}^{-1} r(V^k y) = \sum_{k=0}^\infty r(V^ky) = \infty$, and where $r_n(y)
= \sum_{k=0}^{n-1} r(V^k y)$ if $n\geq 1$, $r_n(y) = 0$ if $n=0$, and $r_n(y) =
\sum_{k=n}^{-1} r(V^k y)$ if $n\leq -1$. For additional intuition and information on
special flows, see Chapter 11 in Cornfeld, Fomin and Sinai
\cite{cornfeld:fomin:sinai:1982}, or Appendix A in Pipiras and Taqqu
\cite{pipiras:taqqu:2003cy}.

\begin{theorem}\label{t:existence-strict-semi-add}
Suppose that $\{F_t\}_{t\in\bbR}$ is an almost semi-additive functional related to a
cocycle $\{A_t\}_{t\in\bbR}$ and to a measurable, nonsingular flow
$\{\phi_t\}_{t\in\bbR}$ on a standard Lebesgue space $(X,{\cal X},\mu)$. Then,
$\{F_t\}_{t\in\bbR}$ has a version which is a semi-additive functional related to the
cocycle $\{A_t\}_{t\in\bbR}$.
\end{theorem}

\begin{proof}
We will extend the proof of Proposition 3.1 in Pipiras and Taqqu
\cite{pipiras:taqqu:2002d} and also use Proposition 1.2 in Kubo \cite{kubo:1970}. By
using Remark 3.1 in Kubo \cite{kubo:1969}, it is enough to prove the proposition in the
following two cases: {\it Case 1:} the flow $\{\phi_t\}_{t\in\bbR}$ is an identity flow,
that is, $\phi_t(x) = x$ for all $t\in\bbR, x\in X$, and {\it Case 2:} the flow
$\{\phi_t\}_{t\in\bbR}$ is a special flow as described above with the function $r$
satisfying $r(y)\geq \theta$ for some fixed $\theta>0$.

\medskip
{\it Case 1:} If the flow $\{\phi_t\}_{t\in\bbR}$ is the identity, then the cocycle $A_t$
in (\ref{e:coc-add}) must satisfy $A_t(x) = 1$ for all $t\in\bbR,x\in X$ (see, for
example, Lemma 3.2 in Pipiras and Taqqu \cite{pipiras:taqqu:2002s}). Relation
(\ref{e:semi-additive-cocycle}) becomes
$$
F_{t_1+t_2}(x) = F_{t_1}(x) + F_{t_2}(\phi_{t_1}(x))\ \ \mbox{a.e.}\ \mu(dx),
$$
for all $t_1,t_2\in\bbR$, which shows that the almost semi-additive functional
$\{F_t\}_{t\in\bbR}$ is also an almost cocycle taking values in $\bbR$ but with addition
as a group operation (it is of the form (\ref{e:coc-add}) but with a sum instead of a
product). Theorem B.9 in Zimmer \cite{zimmer:1984} implies that $\{F_t\}_{t\in\bbR}$ has
a version which is a cocycle and hence, in our terminology, a semi-additive functional
related to the cocycle $\{A_t\}_{t\in\bbR}$.

\medskip
{\it Case 2:} Suppose that $\{\phi_t\}_{t\in\bbR}$ is a special flow on $Y\times
[0,r(\cdot))$ as defined above and hence satisfies (\ref{e:special-representation}). For
notational convenience, we shall write $F(t,(y,u))$ instead of $F_t(y,u)$. Since
$F(t,\cdot)$ is a semi-additive functional related to a cocycle $A_t$, we have that, for
any $s,t\in\bbR$,
$$
F(s+t,(y,u)) = F(s,(y,u)) + A_s(y,u)\, F(t,(V^n y,u+s-r_n(y)))
$$
a.e.\ for $(y,u)$ such that $r_n(y)\leq u+s < r_{n+1}(y)$. We can choose $u = u_0\in
(0,\theta)$ and use the Fubini's theorem to conclude that
\begin{equation}\label{e:main}
F(s+t,(y,u_0)) = F(s,(y,u_0)) + A_s(y,u_0)\, F(t,(V^n y,u_0+s-r_n(y)))
\end{equation}
a.e.\ for $(s,t,y)$ such that $r_n(y)\leq u+s < r_{n+1}(y)$. Setting $n=0$ and $s+u_0=u$
in (\ref{e:main}), we have
\begin{equation}\label{e:main2}
F(t,(y,u)) = (A_{u-u_0}(y,u_0))^{-1} F(t+u-u_0,(y,u_0)) - (A_{u-u_0}(y,u_0))^{-1}
F(u-u_0,(y,u_0))
\end{equation}
a.e.\ for $(t,y,u)$ such that $0\leq u<r(y)$. We shall find expressions for the two
$F$-terms on the right-hand side of (\ref{e:main2}). By using (\ref{e:main}) with ``$s$''
and ``$t$'' indicated by the horizontal braces below, we have
$$
F(t+u-u_0,(y,u_0)) = F(\underbrace{r_n(y)-v}_{s} +
\underbrace{t+u-u_0-r_n(y)+v}_{t},(y,u_0)) = F(r_n(y)-v,(y,u_0))
$$
$$
+\, A_{r_n(y)-v}(y,u_0) F(t+u-u_0-r_n(y)+v,(V^m y,u_0+r_n(y)-v-r_m(y)))
$$
a.e.\ for $(t,y,u,v)$ such that
\begin{equation}\label{e:ineq-mn}
r_m(y) \leq u_0 + r_n(y) - v< r_{m+1}(y).
\end{equation}
We can take $v=v_0\in (0,u_0)$ for which the relation above holds a.e.\ for $(t,y,u)$ so
that we have
\begin{equation}\label{e:uv}
 0<v_0<u_0<\theta<r(y).
\end{equation}
The inequality (\ref{e:ineq-mn}) implies that, for such $v_0$, we have $m=n$ and hence
$$
F(t+u-u_0,(y,u_0)) = F(r_n(y)-v_0,(y,u_0))
$$
\begin{equation}\label{e:aux1}
+ A_{r_n(y)-v_0}(y,u_0) F(t+u-u_0-r_n(y)+v_0,(V^n y,u_0-v_0))
\end{equation}
a.e.\ for $(t,y,u)$. For the second term in (\ref{e:main2}), observe that, by
(\ref{e:main}) with $n=0$,
$$
F(\underbrace{u-u_0}_{t}+\underbrace{v_0}_{s},(y,\underbrace{u_0-v_0}_{u_0})) =
F(v_0,(y,u_0-v_0)) + A_{v_0}(y,u_0-v_0)\, F(u-u_0,(y,u_0))
$$
and hence
\begin{equation}\label{e:aux2}
F(u-u_0,(y,u_0)) = (A_{v_0}(y,u_0-v_0))^{-1} \Big( F(u-u_0+v_0,(y,u_0-v_0)) -
F(v_0,(y,u_0-v_0))\Big)
\end{equation}
a.e.\ for $(y,u)$. Substituting (\ref{e:aux1}) and (\ref{e:aux2}) into (\ref{e:main2}),
we obtain that
\begin{eqnarray}
  F(t,(y,u)) &=& (A_{u-u_0}(y,u_0))^{-1} A_{r_n(y)-v_0}(y,u_0) F(t+u-u_0-r_n(y)+v_0,(V^n y,u_0-v_0))
  \nonumber \\
   & & -\ (A_{u-u_0}(y,u_0))^{-1}  (A_{v_0}(y,u_0-v_0))^{-1}  F(u-u_0+v_0,(y,u_0-v_0)) \nonumber \\
   & & + \ (A_{u-u_0}(y,u_0))^{-1}  F(r_n(y)-v_0,(y,u_0)) \nonumber \\
   & & +\  (A_{u-u_0}(y,u_0))^{-1}  (A_{v_0}(y,u_0-v_0))^{-1} F(v_0,(y,u_0-v_0)) \label{e:aux3}\\
   &=: & F_1(t,(y,u)) + F_2(t,(y,u)) \nonumber
\end{eqnarray}
a.e.\ $(t,y,u)$, where $F_1$ and $F_2$ consist, respectively, of the first and last two
terms  in the sum (\ref{e:aux3}). Since (\ref{e:aux1}) holds for all $n\in\bbZ$, observe
that the relation (\ref{e:aux3}) does not depend on $n\in \bbZ$. We may therefore define
$F_1(t,(y,u))$ and $F_2(t,(y,u))$ as above for $r_n(y)\leq t+u<r_{n+1}(y)$, $n\in\bbZ$.

By using the cocycle equation (\ref{e:coc-add}),
$$
(A_{u-u_0}(y,u_0))^{-1} A_{r_n(y)-v_0}(y,u_0) A_{t+u-r_n(y)-u_0}(V^ny,u_0)
A_{v_0}(V^ny,u_0-v_0)
$$
$$
= (A_{u-u_0}(y,u_0))^{-1} A_{r_n(y)}(y,u_0) A_{t+u-r_n(y)-u_0}(V^ny,u_0)
$$
$$
=(A_{u-u_0}(y,u_0))^{-1} A_{t+u-u_0}(y,u_0) = A_t(y,u),
$$
since $(y,u)\in Y\times [0,r(\cdot))$. Then, we have
\begin{equation}\label{e:F_1}
    F_1(t,(y,u)) = A_{t}(y,u)\, F^*(\phi_t(y,u)) - F^*(y,u),
\end{equation}
where
$$
F^*(y,u) = (A_{u-u_0}(y,u_0))^{-1}  (A_{v_0}(y,u_0-v_0))^{-1}  F(u-u_0+v_0,(y,u_0-v_0)).
$$
Example \ref{ex:semi-additive-cocycle-common} shows that $F_1(t,\cdot)$, given by
(\ref{e:F_1}), is a semi-additive functional related to the cocycle $A_t$. We will now
show that $F_2(t,\cdot)$ can be modified to a semi-additive functional related to the
cocycle $A_t$.

It follows from the definition of the special flow that $\phi_t(y,u) = (y,t+u)$ for
$0\leq u+t <r(y)$, and thus $\phi_t(y,u) = (y,0)$ and $\phi_{v_0}(y,u_0-v_0) = (y,u_0)$
for all $0\leq u_o<r(y)$. Using the cocycle relation (\ref{e:coc-add}), we get
$A_{u-u_0}(y,u_0) = A_{-u_0}(y,u_0)A_u(y,0)$ and $A_{v_0-u_0}(y,u_0-v_0) =
A_{v_0}(y,u_0-v_0)A_{-u_0}(y,u_0)$. Hence,
$$
F_2(t,(y,u)) = (A_u(y,0))^{-1} \Big\{ (A_{-u_0}(y,u_0))^{-1} F(r_n(y)-v_0,(y,u_0))
$$
\begin{equation}\label{e:F_2}
+ (A_{v_0-u_0}(y,u_0-v_0))^{-1} F(v_0,(y,u_0-v_0))\Big\} =: (A_u(y,0))^{-1} G_n(y).
\end{equation}
By Lemma \ref{l:G_n} below,
\begin{equation}\label{e:rec}
G_{n+m}(y) = G_n(y) + A_{r_n(y)}(y,0)\, G_m(V^ny)\quad \mbox{a.e.}\ \mbox{for}\ y.
\end{equation}
It follows that
$$
G_n(y) = \widetilde G_n(y) \ \mbox{a.e.}\ \mbox{for}\ y,
$$
where
\begin{equation}\label{e:G_n}
\widetilde G_n(y) = \sum_{k\in[0,n)}A_{r_k(y)}(y,0)\, G_1(V^k y)
\end{equation}
and $[0,n) = [n,0)$ for $n<0$. This can be checked by using
$$
A_{r_{n+k}(y)}(y,0) = A_{r_n(y)}(y,0) A_{r_k(y)}(\phi_{r_n(y)}(y,0)) =
A_{r_n(y)}(y,0)A_{r_k(y)}(V^ny,0)
$$
to verify that $\widetilde G_n$ satisfies (\ref{e:rec}). By Lemma \ref{l:G_n-cocycle}
below, the function
\begin{equation}\label{e:tilde-F_2}
\widetilde F_2(t,(y,u)) = (A_u(y,0))^{-1} \widetilde G_n(y),
\end{equation}
for $r_n(y)\leq u+t< r_{n+1}(y)$, is a semi-additive functional related to the cocycle
$A_t$.

Since both $F_1$ and $\widetilde F_2$ are semi-additive functionals related to the
cocycle $A_t$, so is the sum
$$
\widetilde F(t,(y,u)) = F_1(t,(y,u)) + \widetilde F_2(t,(y,u))
$$
and, since $F_2(t,(y,u)) = \widetilde F_2(t,(y,u))$ a.e.\ $(t,y,u)$, we have
\begin{equation}\label{e:F=tilde-F}
F(t,(y,u)) = \widetilde F(t,(y,u))
\end{equation}
a.e.\ $(t,y,u)$.

To show that $\{\widetilde F(t,\cdot)\}_{t\in\bbR}$ is a version of
$\{F(t,\cdot)\}_{t\in\bbR}$, it is enough to show that (\ref{e:F=tilde-F}) holds also for
all $t\in\bbR$, a.e.\ $(y,u)$. The argument is standard. Set $\delta(t,(y,u)) =
F(t,(y,u)) - \widetilde F(t,(y,u))$, $\Omega_t = \{(y,u): \delta(t,(y,u)) = 0\}$ and also
$\Omega_{s,t} = \{(y,u): \delta(s+t,(y,u)) = \delta(s,(y,u)) + A_s(y,u)
\delta(t,\phi_s(y,u))\}$. Denoting the Lebesgue measure on $\bbR$ by $\bbL$, we have
$(\tau\otimes\bbL)(\Omega_{s,t}^c) = 0$ for all $s,t\in\bbR$ but only
$(\tau\otimes\bbL)(\Omega_t^c) = 0$ a.e.\ for $t\in\bbR$. However, if $r\in\bbR$, then
there are $s,t\in\bbR$ such that $s+t=r$ and $(\tau\otimes\bbL)(\Omega_{s}^c) =
(\tau\otimes\bbL)(\Omega_{t}^c) = 0$. We also have
$(\tau\otimes\bbL)((\phi_{-s}\Omega_{t})^c) = 0$ since $\phi_{s}$ is one-to-one and onto.
Then, $(\tau\otimes\bbL)((\Omega_{s,t}\cap \Omega_s \cap \phi_{-s}\Omega_t)^c) = 0$ and,
for $(y,u)\in \Omega_{s,t}\cap \Omega_s \cap \phi_{-s}\Omega_t$, we have $\delta(r,(y,u))
= \delta(s,(y,u)) + A_s(y,u) \delta(t,\phi_s(y,u)) =0$. This shows that
$(\tau\otimes\bbL)\{F(r,(y,u)) \neq \widetilde F(r,(y,u)\} =0$ for any $r\in\bbR$, that
is, $\{\widetilde F(t,\cdot)\}_{t\in\bbR}$ is a version of $\{F(t,\cdot)\}_{t\in\bbR}$. \
\ $\Box$
\end{proof}

\medskip

The theorem can be readily expressed in terms of multiplicative flows by changing time
$t\in\bbR$ into $c = e^t$, $c>0$.

\begin{corollary}\label{c:existence-strict-semi-add-m}
Suppose that $\{J_c\}_{c>0}$ is an almost semi-additive functional related to a cocycle
$\{B_c\}_{c>0}$ satisfying (\ref{e:semi-additive-cocycle-m}) and to a measurable,
nonsingular multiplicative flow $\{\psi_c\}_{c>0}$ on a standard Lebesgue space $(X,{\cal
X},\mu)$. Then, $\{J_c\}_{c>0}$ has a version which is a semi-additive functional related
to the cocycle $\{B_c\}_{c>0}$.
\end{corollary}

\medskip
The next two corollaries, the first one for additive flows and the second for
multiplicative flows, provide insight into the structure of semi-additive functionals
related to a cocycle. Corollary \ref{c:solutions-m} will be used in Propositions
\ref{p:1-semi-add-cyclic} and \ref{p:2-semi-add-cyclic} below to deduce the forms of the
{\it 1}- and {\it 2}-semi-additive functionals corresponding to a cyclic flow.

\begin{corollary}\label{c:solutions}
If the flow $\{\phi_t\}_{t\in\bbR}$ is given by its special representation
(\ref{e:special-representation}) on a space $Y\times [0,r(\cdot))$, with the maps $V$ and
$r$, and $\{F_t\}_{t\in\bbR}$ is a semi-additive functional related to a cocycle
$\{A_t\}_{t\in\bbR}$, then
\begin{equation}\label{e:solution-semi-additive}
F_t(y,u) = F^{(1)}_t(y,u) + F^{(2)}_t(y,u),
\end{equation}
where
\begin{eqnarray}
  F^{(1)}_t(y,u) &=& A_t(y,u)\, F(\phi_t(y,u)) - F(y,u),  \label{e:solution-1}\\
  F^{(2)}_t(y,u) &=& (A_u(y,0))^{-1} \sum_{k\in[0,n)} A_{r_k(y)}(y,0) F_1(V^ky),
  \label{e:solution-2}
\end{eqnarray}
for $r_n(y)\leq t+u< r_{n+1}(y)$, $F,F_1$ are some functions and $[0,n) = [n,0)$ for
$n<0$. Moreover, each of the functions $F^{(1)}(t,\cdot)$ and $F^{(2)}(t,\cdot)$ is a
semi-additive functional related to the cocycle $\{A_t\}_{t\in\bbR}$.
\end{corollary}

\begin{proof}
The corollary follows from the proof of Theorem \ref{t:existence-strict-semi-add} by
replacing ``a.e.'' by ``for all'' conditions. See, in particular, (\ref{e:aux3}),
(\ref{e:F_1}) and (\ref{e:F_2}) together with (\ref{e:G_n}). \ \ $\Box$
\end{proof}

\medskip
Corollary \ref{c:solutions} involves an additive flow $\{\phi_t\}_{t\in\bbR}$. The next
corollary formulates the result for a multiplicative flow $\{\psi_c\}_{c>0}$.

\begin{corollary}\label{c:solutions-m}
If $\{J_c\}_{c>0}$ is a semi-additive functional related to the cocycle $\{B_c\}_{c>0}$
and a {\em multiplicative} flow $\{\psi_c\}_{c>0}$, given by its special representation
$\psi_c(y,u) = (V^ny,u+\ln c-r_n(y))$, then
\begin{equation}\label{e:solution-semi-additive-m}
J_c(y,u) = J^{(1)}_c(y,u) + J^{(2)}_c(y,u),
\end{equation}
where
\begin{eqnarray}
  J^{(1)}_c(y,u) &=& B_c(y,u)\, J(\psi_c(y,u)) - J(y,u),  \label{e:solution-1-m}\\
  J^{(2)}_c(y,u) &=& (B_{e^u}(y,0))^{-1} \sum_{k\in[0,n)} B_{e^{r_k(y)}}(y,0) J_1(V^ky),
  \label{e:solution-2-m}
\end{eqnarray}
for $r_n(y)\leq \ln c + u< r_{n+1}(y)$, and $J,J_1$ are some functions.
\end{corollary}

\begin{proof}
This result follows by observing that $F_t(y,u):= J_{e^t}(y,u)$ is a semi-additive
functional related to the cocycle $A_t(y,u) = B_{e^t}(y,u)$ and the additive flow
$\phi_t(y,u) = \psi_{e^t}(y,u)$, applying Corollary \ref{c:solutions} to $F_t(y,u)$ and
then translating the result back in terms of the map $J_c(y,u)$. \ \ $\Box$
\end{proof}

\medskip
The following two auxiliary lemmas were used in the proof of Theorem
\ref{t:existence-strict-semi-add}.

\begin{lemma}\label{l:G_n}
Let $G_n(y)$ be defined by (\ref{e:F_2}). Then, for $n,m\in\bbZ$,
\begin{equation}\label{e:G_n+m}
G_{n+m}(y) = G_n(y) + A_{r_n(y)}(y,0) G_m(V^ny)\ \mbox{a.e.}\ y.
\end{equation}
\end{lemma}

\begin{proof}
Recall from (\ref{e:F_2}) that
\begin{equation}\label{e:G_n-again}
G_n(y) = (A_{-u_0}(y,u_0))^{-1} F(r_n(y)-v_0,(y,u_0)) + (A_{v_0-u_0}(y,u_0-v_0))^{-1}
F(v_0,(y,u_0-v_0))
\end{equation}
where $u_0$ and $v_0$ satisfy (\ref{e:uv}). To show (\ref{e:G_n+m}), we shall use the
relation
\begin{equation}\label{e:r_n+m}
r_{n+m}(y) = r_n(y) + r_m(V^ny),
\end{equation}
which is easy to verify by using the definition of $r_n(y)$. Observe that, for
$n,m\in\bbZ$, by using (\ref{e:r_n+m}) and (\ref{e:main}),
$$
F(r_{n+m}(y)-v_0,(y,u_0)) =  F(\underbrace{r_{n}(y)-v_0}_{s} +
\underbrace{r_m(V^ny)}_{t},(y,u_0))
$$
$$
= F(r_{n}(y)-v_0,(y,u_0)) + A_{r_{n}(y)-v_0}(y,u_0)\, F(r_m(V^ny),(V^ny,u_0-v_0))
$$
$$
= F(r_{n}(y)-v_0,(y,u_0)) + A_{r_{n}(y)-v_0}(y,u_0)\, F(\underbrace{v_0}_{s} +
\underbrace{r_m(V^ny)-v_0}_{t},(V^ny,u_0-v_0))
$$
$$
= F(r_{n}(y)-v_0,(y,u_0)) + A_{r_{n}(y)-v_0}(y,u_0)\, F(v_0,(V^ny,u_0-v_0))
$$
\begin{equation}\label{e:F_r_n+m}
+ A_{r_{n}(y)-v_0}(y,u_0)A_{v_0}(V^ny,u_0-v_0)\, F(r_m(V^ny)-v_0,(V^ny,u_0)).
\end{equation}
To compute $G_{n+m}$ we use (\ref{e:G_n-again}) and substitute (\ref{e:F_r_n+m}) for
$F(r_{n+m}(y)-v_0,(y,u_0))$. This yields
$$
G_{n+m}(y) = (A_{-u_0}(y,u_0))^{-1} F(r_n(y)-v_0,(y,u_0)) + (A_{v_0-u_0}(y,u_0-v_0))^{-1}
F(v_0,(y,u_0-v_0))
$$
$$
+ (A_{-u_0}(y,u_0))^{-1} A_{r_{n}(y)-v_0}(y,u_0)A_{v_0}(V^ny,u_0-v_0)\,
F(r_m(V^ny)-v_0,(V^ny,u_0))
$$
$$
+(A_{-u_0}(y,u_0))^{-1} A_{r_{n}(y)-v_0}(y,u_0)\, F(v_0,(V^ny,u_0-v_0))
$$
$$
= G_n(y) + (A_{-u_0}(y,u_0))^{-1}
A_{r_{n}(y)-v_0}(y,u_0)A_{v_0}(V^ny,u_0-v_0)A_{-u_0}(V^ny,u_0) \cdot
$$
$$
\cdot (A_{-u_0}(V^ny,u_0))^{-1} F(r_m(V^ny)-v_0,(V^ny,u_0))
$$
$$
+ (A_{-u_0}(y,u_0))^{-1} A_{r_{n}(y)-v_0}(y,u_0) A_{v_0-u_0}(V^ny,u_0-v_0)\cdot
$$
$$
\cdot (A_{v_0-u_0}(V^ny,u_0-v_0))^{-1} F(v_0,(V^ny,u_0-v_0))
$$
\begin{equation}\label{e:Gend}
= G_n(y) + (A_{-u_0}(y,u_0))^{-1} A_{r_{n}(y)-v_0}(y,u_0) A_{v_0-u_0}(V^ny,u_0-v_0)\,
G_{m}(V^n y)
\end{equation}
a.e.\ $y$, where to obtain the last identity we used the relation
$A_{v_0}(V^ny,u_0-v_0)A_{-u_0}(V^ny,u_0) = A_{v_0-u_0}(V^ny,u_0-v_0)$. Since
$A_{-u_0}(y,u_0) A_{u_0}(y,0) = A_{u_0-u_0}(y,0) = A_0(y,0) = 1$ (see (\ref{e:a-zero})),
we get $(A_{-u_0}(y,u_0))^{-1} = A_{u_0}(y,0)$ and, since $A_{u_0}(y,0)
A_{r_{n}(y)-v_0}(y,u_0) A_{v_0-u_0}(V^ny,u_0-v_0) = A_{r_{n}(y)-v_0 + u_0}(y,0)
A_{v_0-u_0}(V^ny,u_0-v_0) = A_{r_n(y)}(y,0)$, we see that (\ref{e:Gend}) reduces to
(\ref{e:G_n+m}).\ \ $\Box$
\end{proof}

\begin{lemma}\label{l:G_n-cocycle}
If $F_t(y,u)$ is defined by the right-hand side of (\ref{e:tilde-F_2}), together with
(\ref{e:G_n}), then $F_t(y,u)$ is a semi-additive functional related to the cocycle
$\{A_t\}_{t\in\bbR}$.
\end{lemma}

\begin{proof}
Fix $s,t\in\bbR$ and $(y,u)\in Y\times [0,r(\cdot))$. We are interested in $F_s(y,u)
=(A_{u}(y,0))^{-1} \widetilde G_n(y)$ where $n$ is such that $r_n(y) \leq s+u<
r_{n+1}(y)$ and in $F_{s+t}(y,u) =(A_{u}(y,0))^{-1} \widetilde G_{n+m}(y)$ where, in
addition, $m$ is such that $r_{n+m}(y) \leq t+s+u< r_{n+m+1}(y)$, an inequality
equivalent to $r_m(V^ny) \leq t+s+u-r_n(y)< r_{m+1}(V^ny)$ by (\ref{e:r_n+m}). Then, by
using the fact that $\widetilde G_n$ satisfies (\ref{e:G_n+m}), we obtain that
$$
A_s(y,u) F_t(\phi_s(y,u)) = A_s(y,u) F_t(V^ny,s+u-r_n(y)) = A_s(y,u)
(A_{s+u-r_n(y)}(V^ny,0))^{-1} \widetilde G_m(V^ny)
$$
$$
= A_s(y,u) (A_{s+u-r_n(y)}(V^ny,0))^{-1} (A_{r_n(y)}(y,0))^{-1} \Big( \widetilde
G_{n+m}(y) - \widetilde G_n(y)\Big)
$$
$$
= A_s(y,u) (A_{s+u}(y,0))^{-1} \Big( \widetilde G_{n+m}(y) - \widetilde G_n(y)\Big) =
(A_{u}(y,0))^{-1} \Big( \widetilde G_{n+m}(y) - \widetilde G_n(y)\Big)
$$
$$
= F_{s+t}(y,u) - F_s(y,u).
$$
This concludes the proof. \ \ $\Box$
\end{proof}


\section{{\it 2}-semi-additive functionals}
\label{s:new-functional}

In this section, we apply Corollary \ref{c:existence-strict-semi-add-m} to show that the
function $j_c(x)$ in the relation (\ref{e:generated}) can be chosen, without loss of
generality, as a {\it 2}-semi-additive functional satisfying (\ref{e:II-semi-add}).
Observe first that because of the properties of the Radon-Nikodym derivatives, one has
for all $c_1,c_2>0$,
\begin{eqnarray}
  \frac{d(\mu\circ \psi_{c_1c_2})}{d\mu}(x) &=& \frac{d(\mu\circ \psi_{c_1})}{d\mu}(x)
\frac{d(\mu\circ \psi_{c_2} \circ \psi_{c_1})}{d(\mu \circ \psi_{c_1})}(x)\nonumber \\
  &=& \frac{d(\mu\circ \psi_{c_1})}{d\mu}(x) \frac{d(\mu\circ \psi_{c_2})}{d\mu}
(\psi_{c_1}(x)) \label{e:radon-nikodym-cocycle}
\end{eqnarray}
a.e.\ $\mu(dx)$, which is the relation (\ref{e:cocycle}) defining a cocycle but valid
only a.e.\ $\mu(dx)$ and not for all $x\in X$. We start with the following lemma which
shows that the collection of the Radon-Nikodym derivatives $\{ d(\mu\circ
\psi_c)/d\mu\}_{c>0}$ has a version which is a cocycle for all $x\in X$.

\begin{lemma}\label{l:radon-nikodym-cocycle}
Suppose that the relations (\ref{e:generated}) and (\ref{e:full-support}) hold. Then, the
Radon-Nikodym derivatives $\{ d(\mu\circ \psi_c)/d\mu\}_{c>0}$ have a version which is

1) jointly measurable in $(c,x)$,

2) a cocycle mapping $(0,\infty)\times X \to (0,\infty)$,

3) a Radon-Nikodym derivative $d(\mu\circ \psi_c)/d\mu$ for all $c>0$.
\end{lemma}

\begin{proof}
We first show that the collection $\{ d(\mu\circ \psi_c)/d\mu\}_{c>0}$ has a jointly
measurable version. By using the notation (\ref{e:G_t}), relation (\ref{e:generated})
implies that, for any $t\in\bbR$ and $c>0$,
\begin{equation}\label{e:generated-G_t}
G_{ct}(x,cu) = c^{H-1/\alpha} b_c(x)  \left\{\frac{d(\mu\circ \psi_{c})}{d\mu}(x)
\right\}^{1/\alpha} G_t(\psi_c(x), u + g_c(x))
\end{equation}
a.e.\ $\mu(dx)du$.

If $\mbox{supp} \{G_t(x,u)\}= X \times \bbR$ a.e.\ $\mu(dx)du$ for some fixed $t\in
\bbR$, we have for $c>0$,
$$
\frac{d(\mu\circ \psi_{c})}{d\mu}(x) =  \left\{ \frac{c^{1/\alpha -  H}G_{ct}(x,cu)}{
b_c(x) G_t(\psi_c(x), u + g_c(x))} \right\}^\alpha
$$
a.e.\ $\mu(dx)du$. Hence, since the right-hand side of the expression above is jointly
measurable, we may conclude that $\{ d(\mu\circ \psi_c)/d\mu\}_{c>0}$ has a jointly
measurable version.

Consider now the general case when $\mbox{supp} \{G_t(x,u)\}= X \times \bbR$ a.e.\
$\mu(dx)du$ may possibly not hold for any $t\in \bbR$. Let ${\cal G} =
\mbox{Sp}\{G_t,t\in\bbR\}$ be the linear span of $G_t$, $t\in\bbR$, and $\overline {\cal
G}$ be the closure of ${\cal G}$ in the space $L^\alpha(X\times\bbR,\mu(dx)du)$. Since
$\{G_t,t\in\bbR\}\subset \overline {\cal G}$, the assumption (\ref{e:full-support})
implies that $\mbox{supp}\{\overline {\cal G}\} = X\times \bbR$ a.e.\ $\mu(dx)du$. By
Lemma 3.2 in Hardin \cite{hardin:1981}, there is a function $G^*\in \overline {\cal G}$
such that $\mbox{supp} \{G^*(x,u)\}= X \times \bbR$ a.e.\ $\mu(dx)du$. Since $G^*\in
\overline {\cal G}$, there are functions $G^{(n)}(x,u) = \sum_i a_{ni} G_{t_{ni}}(x,u)\in
{\cal G}$, $n\geq 1$, $a_{ni},t_{ni}\in\bbR$, such that $G^{(n)}(x,u) \to G^*(x,u)$ a.e.\
$\mu(dx)du$.

Let also $G^{(n)}_c(x,u) = \sum_i a_{ni} G_{ct_{ni}}(x,cu)$. Relation
(\ref{e:generated-G_t}) implies that, for any $c>0$ and $n\geq 1$,
\begin{equation}\label{e:generated-G^n}
G^{(n)}_c(x,u) = c^{H-1/\alpha} b_c(x) \left\{\frac{d(\mu\circ \psi_{c})}{d\mu}(x)
\right\}^{1/\alpha}  G^{(n)}(\psi_c(x), u + g_c(x))
\end{equation}
a.e.\ $\mu(dx)du$. For any $c>0$, the right-hand side of (\ref{e:generated-G^n})
converges to
$$
c^{H-1/\alpha} b_c(x) \left\{\frac{d(\mu\circ \psi_{c})}{d\mu}(x) \right\}^{1/\alpha}
G^*(\psi_c(x), u + g_c(x))
$$
a.e.\ $\mu(dx)du$, as $n\to\infty$. Since the right-hand side of (\ref{e:generated-G^n})
converges, the left-hand side of (\ref{e:generated-G^n}) converges to some function
$G^*_c(x,u)$. Hence, for any $c>0$,
$$
G^*_c(x,u) = c^{H-1/\alpha} b_c(x) \left\{\frac{d(\mu\circ \psi_{c})}{d\mu}(x)
\right\}^{1/\alpha}  G^*(\psi_c(x), u + g_c(x))
$$
or, since $\mbox{supp} \{G^*\}= X \times \bbR$ a.e.,
\begin{equation}\label{e:generated-G^*}
\frac{d(\mu\circ \psi_{c})}{d\mu}(x) =  \left\{ \frac{c^{1/\alpha -  H}G_{c}^*(x,u)}{
b_c(x) G^*(\psi_c(x), u + g_c(x))} \right\}^\alpha
\end{equation}
a.e.\ $\mu(dx)du$. Observe that $G^*_c(x,u)$ is jointly measurable in $(c,x,u)$ because
it is the a.e.\ limit of functions jointly measurable in $(c,x,u)$. Since $G^*(x,u)$ is
measurable in $(x,u)$, $\psi_c(x),g_c(x)$ and $b_c(x)$ are measurable in $(c,x)$, the
function $c^{H-1/\alpha} b_c(x)^{1/\alpha} G^*(\psi_c(x), u + g_c(x))$ is measurable in
$(c,x,u)$. Hence, the right-hand side of (\ref{e:generated-G^*}) is jointly measurable
which is to say that $\{d(\mu\circ \psi_c)/d\mu\}_{c>0}$ has a jointly measurable
version.

Suppose then, without loss of generality, that $d(\mu\circ \psi_c)/d\mu(x)$ is jointly
measurable. We still need to show that $\{d(\mu\circ \psi_c)/d\mu\}_{c>0}$ has a version
which is a cocycle. Since the flow $\{\psi_c\}_{c>0}$ is nonsingular, the measures $\mu
\circ \psi_c$ and $\mu$ are equivalent and hence we may suppose that $(d(\mu\circ
\psi_c)/d\mu)(x):(0, \infty) \times X \to \bbR \setminus \{0\}.$ By
(\ref{e:radon-nikodym-cocycle}), $\{d(\mu\circ \psi_c)/d\mu\}_{c>0}$ is an {\it almost
cocycle} for the flow $\{\psi_c\}_{c>0}$ where ``almost'' refers to the fact that the
relation (\ref{e:radon-nikodym-cocycle}) holds a.e.\ $\mu (dx)$ for $c_1,c_2>0$, in
contrast to (\ref{e:cocycle}) which holds for all $x \in X$ and $c_1,c_2>0$. By Theorem
B.9 in Zimmer \cite{zimmer:1984} (see also Theorem A.1 in Kolody{\' n}ski and Rosi{\'
n}ski \cite{kolodynski:rosinski:2002}) and since $(d(\mu\circ \psi_c)/d\mu)(x)$ is
measurable in $(c,x)$, $\{d(\mu\circ \psi_c)/d\mu\}_{c>0}$ has a version which is a
cocycle for the flow $\{ \psi_c \}_{c>0}$ taking values in $(0,\infty)$. Property $3)$
follows from the definition of ``version''.
 \ \ $\Box$
\end{proof}

\medskip

\noindent {\bf Remark.} The version specified in Lemma \ref{l:radon-nikodym-cocycle}
which satisfies Conditions 1), 2) and 3) is not unique. Suppose, for instance, that $X
=\bbR^2 = \{(x_1,x_2)\}$, $\mu(dx) = dx_1dx_2$ and $\psi_{c}(x_1,x_2) = (x_1,x_2+\ln c)$.
Then,
$$
\frac{d(\mu\circ \psi_{c})}{d\mu}(x_1,x_2) \equiv 1
$$
is a version of the Radon-Nikodym derivatives satisfying Conditions 1), 2) and 3). On the
other hand, let $b:\bbR \mapsto (0,\infty)$ be an arbitrary function and $x_1^*\in\bbR$
be fixed. Then,
$$
\frac{d(\mu\circ \psi_{c})}{d\mu}(x_1,x_2) = \left\{
\begin{array}{cc}
 1, & x_1\neq x_1^*, \\
 \frac{b(x_2 + \ln c)}{b(x_2)}, & x_1 = x_1^*,
 \end{array}\right.
$$
is also a version of the Radon-Nikodym derivatives satisfying Conditions 1), 2) and 3).
Indeed, it is jointly measurable and also, for fixed $c>0$, it is still a Radon-Nikodym
derivative since it was modified on the set $\{(x_1,x_2):x_1 = x_1^*\}$ of a
$\mu$-measure zero. It satisfies a cocycle equation for all $x\in \bbR^2$, $c_1,c_2>0$,
because it does so on the disjoint subsets $\{(x_1,x_2):x_1\neq x_1^*\}$ and
$\{(x_1,x_2):x_1= x_1^*\}$ of $\bbR^2$ which are invariant under the flow. Observe that
the two versions of the Radon-Nikodym derivatives above are different when $b\not\equiv
1$.

\medskip
We can now give a precise definition of {\it 2}-semi-additive functional.

\begin{definition}\label{d:1-2-semi-additive} A measurable
function $j_c(x):(0, \infty) \times X \to \bbR$ is a {\it 2-semi-additive functional} for
a flow $\{\psi_c\}_{c>0}$ and a cocycle $\{b_c\}_{c>0}$ if the relation
(\ref{e:II-semi-add}) holds, where $d(\mu \circ \psi_c)/d \mu$ satisfies Conditions 1),
2) and 3) of Lemma \ref{l:radon-nikodym-cocycle}. A semi-additive functional
$\{g_c\}_{c>0}$ satisfying (\ref{e:I-semi-add}) is called a {\it 1-semi-additive
functional}.
\end{definition}

In the following examples, we show that after multiplication by a suitable factor, the
{\it 1}-semi-additive functional $\{g_c\}_{c>0}$ in (\ref{e:I-semi-add}) and the {\it
2}-semi-additive functionals $\{j_c\}_{c>0}$ in (\ref{e:II-semi-add}) become
semi-additive functionals related to a cocycle, and we identify these cocycles.

\begin{example}\label{ex:semi-additive-cocycle2}
If $\{g_c\}_{c>0}$ is a {\it 1}-semi-additive functional satisfying (\ref{e:I-semi-add}),
then $J_c(x) = c g_c(x)$ satisfies
$$
J_{c_1c_2}(x) = J_{c_1}(x) + c_1 J_{c_2}(\psi_{c_1}(x)), \quad \mbox{for all} \
c_1,c_2>0, x\in X.
$$
But $B_c(x) = c$ is a cocycle for the flow $\{\psi_c\}_{c>0}$ (it satisfies
(\ref{e:cocycle})). Therefore, $\{J_c\}_{c>0}$ is a semi-additive functional related to
the cocycle $\{B_c\}_{c>0}$.
\end{example}

\begin{example}\label{ex:semi-additive-cocycle1}
If $\{j_c\}_{c>0}$ is a {\it 2}-semi-additive functional satisfying
(\ref{e:II-semi-add}), then $J_c(x)= c^{H-1/\alpha} j_c(x)$ satisfies
$$
J_{c_1c_2}(x) = J_{c_1}(x) + B_{c_1}(x) J_{c_2}(\psi_{c_1}(x)), \quad \mbox{for all} \
c_1,c_2>0, x\in X$$ with
\begin{equation}\label{e:cocycle-2-semi-additive-again}
    B_c(x) = c^{H-1/\alpha} b_c(x)
    \left\{\frac{d(\mu\circ\psi_c)}{d\mu}(x)\right\}^{1/\alpha}.
\end{equation}
Since $\{b_c\}_{c>0}$ is a cocycle taking values in $\{-1,1\}$ and $\{d(\mu\circ
\psi_c)/d\mu\}_{c>0}$ is a cocycle taking values in $\bbR\setminus \{0\}$, it is easy to
check that $\{B_c\}_{c>0}$ is also a cocycle taking values in $\bbR\setminus \{0\}$.
Thus, $\{J_c\}_{c>0}$ is a semi-additive functional related to the cocycle
(\ref{e:cocycle-2-semi-additive-again}).
\end{example}

\noindent {\bf Remarks}
\begin{enumerate}

\item The cocycles $\{B_c\}_{c>0}$ appearing in the preceding examples are associated
with functionals $\{J_c\}_{c>0}$. They should not be confused with the cocycle
$\{b_c\}_{c>0}$ in Relation (\ref{e:generated}).

\item The preceding examples can be used in the following way. Suppose that
$\{g_c\}_{c>0}$ and $\{j_c\}_{c>0}$ are only {\it almost} semi-additive functionals. Then
$\{J_c\}_{c>0}$ would also be almost semi-additive functionals. Since Corollary
\ref{c:solutions-m} applies to $\{J_c\}_{c>0}$, these have a version which is a
semi-additive functional. In view of the expressions relating $J_c$ to $g_c$ and $j_c$,
it follows that $\{g_c\}_{c>0}$ and $\{j_c\}$ have also a version which is a
semi-additive functional. This type of argument is used in the proof of the following
theorem.

\end{enumerate}

\begin{theorem}\label{t:2-semi-additive} Let $\alpha\in (0,2)$ and $H>0$. The
function $j_c(x)$ in relation (\ref{e:generated}) can be taken to be a {\it
2}-semi-additive functional.
\end{theorem}

\begin{proof}
We need first to show that the function $j_c(x)$ in (\ref{e:generated}) is an almost
semi-additive function. Observe that it equals
\begin{equation}\label{e:2-semi-additive-G}
    j_c(x) = c^{-(H-1/\alpha)} G(x,cu) - \widetilde b_c(x)\, G(\psi_c(x),u+g_c(x)),
\end{equation}
a.e.\ $\mu(dx)du$, for any $c>0$, where
\begin{equation}\label{e:cocycle-2-semi-additive}
    \widetilde b_c(x) =  b_c(x)
    \left\{\frac{d(\mu\circ\psi_c)}{d\mu}(x)\right\}^{1/\alpha}.
\end{equation}
By Lemma \ref{l:radon-nikodym-cocycle} above, the Radon-Nikodym derivative
$d(\mu\circ\psi_c)/d\mu$ and hence its $1/\alpha$-power have a version which is a cocycle
taking values in $(0,\infty)$. Since $b_c(x)$ is a cocycle, the product $\widetilde
b_c(x)$ also has a version which is a cocycle taking values in $\bbR\setminus \{0\}$. We
may therefore suppose without loss of generality that $\{\widetilde b_c\}_{c>0}$ is a
cocycle in (\ref{e:2-semi-additive-G}), that is,
\begin{equation}\label{e:bcoc}
\widetilde b_{c_1c_2}(x) = \widetilde b_{c_1}(x)\widetilde b_{c_2}(\psi_{c_1}(x)).
\end{equation}
By using (\ref{e:2-semi-additive-G}), we have for $c_1,c_2>0$,
$$
j_{c_1c_2}(x)  = (c_1 c_2)^{-(H-1/\alpha)} G(x,c_1c_2u) - \widetilde b_{c_1c_2}(x)\,
G(\psi_{c_1c_2}(x),u+g_{c_1c_2}(x))
$$
a.e.\ $\mu(dx)du$. By using (\ref{e:bcoc}) and $g_{c_1c_2}(x) = c_2^{-1}g_{c_1}(x) +
g_{c_2}(\psi_{c_1}(x))$, we conclude that
\begin{eqnarray}
 j_{c_1c_2}(x)  &=& c_2^{-(H-1/\alpha)} \Big( c_1^{-(H-1/\alpha)} G(x,c_1(c_2u)) -
 \widetilde b_{c_1}(x)\, G(\psi_{c_1}(x),c_2u+g_{c_1}(x))
 \Big) \nonumber \\
 & & \hspace{0.5in} + \
 \widetilde b_{c_1}(x) \bigg\{ c_2^{-(H-1/\alpha)} G\Big(\psi_{c_1}(x),c_2(u+c_2^{-1}g_{c_1}(x))\Big)  \nonumber \\
 & & \hspace{1in} -\ \widetilde b_{c_2} (\psi_{c_1} (x)) G\Big(\psi_{c_2}(\psi_{c_1}(x)), u + c_2^{-1}g_{c_1}(x) +
g_{c_2}(\psi_{c_1}(x))\Big)\bigg\} \nonumber \\
 &=& c_2^{-(H-1/\alpha)} j_{c_1}(x) +
 \widetilde b_{c_1}(x) j_{c_2}(\psi_{c_1}(x)) \label{e:2-semi-additive=almost-semi-additive}
\end{eqnarray}
a.e.\ $\mu(dx)$. Hence, $\{j_c\}_{c>0}$ is an almost semi-additive functional.

Multiplying (\ref{e:2-semi-additive=almost-semi-additive}) by $(c_1 c_2)^{H-1/\alpha}$
and setting $J_c(x)=c^{H-1/\alpha} j_c(x)$, $B_c(x)=c^{H-1/\alpha} \widetilde b_c(x)$, we
obtain that
\begin{equation} \label{e:widetilde-j-c1c2}
J_{c_1,c_2}(x)= J_{c_1}(x) + B_{c_1}(x) J_{c_2}(\psi_{c_1}(x)) \quad \mbox{a.e.} \
\mu(dx).
\end{equation}
Since $\{B_c\}_{c>0}$ is also a cocycle for the flow $\{\psi_c\}_{c>0}$, relation
(\ref{e:widetilde-j-c1c2}) shows that $\{J_c\}_{c>0}$ is an almost semi-additive
functional related to the cocycle $\{B_c\}_{c>0}$ in the sense of Definition
\ref{d:semi-additive-cocycle} below. By Corollary \ref{c:existence-strict-semi-add-m},
$\{J_c\}_{c>0}$ has a version which is a semi-additive functional related to the cocycle
$\{B_c\}_{c>0}$. But $j_c(x) = c^{-(H-1/\alpha)} J_c(x)$. Hence, when multiplied by
$c^{-(H-1/\alpha)}$, this version is a {\it 2}-semi-additive functional which is a
version of $\{j_c\}_{c>0}$. \ \ $\Box$
\end{proof}

\medskip
The corresponding result for $g_c$ was proved in Proposition 3.1 of Pipiras and Taqqu
\cite{pipiras:taqqu:2002d}, pages 421-426. (In that proposition, ``semi-additive
functional'' means ``{\it 1}-semi-additive-functional''.) The proof of that Proposition
3.1 also follows from the more general setting of the present paper. In fact, it reduces,
at this stage, to the argument in Example \ref{ex:semi-additive-cocycle2} and the remark
preceding the statement of Theorem \ref{t:2-semi-additive}.

Theorem \ref{t:2-semi-additive} is useful when solving for the kernel function $G$
generated by a given flow. One such case, studied in Proposition 3.1 of Pipiras and Taqqu
\cite{pipiras:taqqu:2003in}, concerns kernel functions related to cyclic flows. {\it
2}-semi-additive functionals could have also been used in Theorem 5.1 of Pipiras and
Taqqu \cite{pipiras:taqqu:2002s}, where kernels related to identity flows of Example
\ref{ex:2-semi-additive} below are considered. In particular, the argument used in
Theorem 5.1 of that paper involving an almost everywhere version of the Cauchy functional
equation would not be necessary anymore. We explain this in greater detail in the remark
following Example \ref{ex:1-semi-additive-ident-dissip} below.


\section{Examples}
\label{s:examples}

In the following examples, we consider {\it 1}- and {\it 2}-semi-additive functionals
$\{j_c\}_{c>0}$ for identity, dissipative and cyclic flows $\{\psi_c\}_{c>0}$, and
related cocycles $\{b_c\}_{c>0}$.

\begin{example}\label{ex:2-semi-additive}
Consider the {\it identity} flow $\{\psi_c\}_{c>0}$ on $(X,\mu)$ such that
\begin{equation}\label{e:identity-flow}
\psi_c(x)=x
\end{equation}
for all $c>0$ and $x\in X$. We can take
\begin{equation}\label{e:jacobian-identity}
\frac{d(\mu\circ\psi_c)}{d\mu}(x) = \frac{d\mu}{d\mu}(x) \equiv 1
\end{equation}
as a cocycle for the identity flow $\{\psi_c\}_{c>0}$. By Lemma 3.2 in Pipiras and Taqqu
\cite{pipiras:taqqu:2002s},
\begin{equation}\label{e:cocycle-identity}
b_c(x)=1
\end{equation}
for the identity flow $\{\psi_c\}_{c>0}$. The {\it 2}-semi-additive functional
$\{j_c\}_{c>0}$ in (\ref{e:II-semi-add}) therefore satisfies
\begin{equation}\label{e:2-semi-additive-identity}
    j_{c_1c_2}(x) = c_2^{-(H-1/\alpha)} j_{c_1}(x) + j_{c_2}(x),
\end{equation}
for all $x\in X$, $c_1,c_2>0$. Relation (\ref{e:2-semi-additive-identity}) is an equation
for the functional $j_c(x)$, which we shall now solve.

If $H\neq 1/\alpha$, by subtracting $ j_{c_1c_2}(x) = c_1^{-(H-1/\alpha)} j_{c_2}(x) +
j_{c_1}(x)$ from (\ref{e:2-semi-additive-identity}), we obtain that
$$
(1-c_2^{-(H-1/\alpha)}) j_{c_1}(x) = (1-c_1^{-(H-1/\alpha)}) j_{c_2}(x).
$$
By fixing $c_2\neq 1$, we conclude that $$ j_c(x) = j(x) (1-c^{-(H-1/\alpha)}), $$ where
$j(x)$ is some function. If $H=1/\alpha$, then
$$
j_{c_1c_2}(x) = j_{c_2}(x) + j_{c_1}(x)
$$
and by using Lemma 1.1.6 in Bingham et al.\ \cite{bingham:goldie:teugels:1987}, we have
$$ j_c(x) = j(x) \ln c,$$ where $j(x)$ is some function.
\end{example}

\begin{example}\label{ex:2-semi-additive-dissipative}
Consider the flow
\begin{equation}\label{e:diss}
\psi_c(y,u) = (y,u+\ln c)
\end{equation}
on the space $(X,\mu) = (Y\times \bbR,\nu(dy)du)$. By Krengel's theorem (see, for
example, Theorem 3.1 in Pipiras and Taqqu \cite{pipiras:taqqu:2002s}), any {\it
dissipative} flow on $(X,\mu)$ is null-isomorphic to a flow $\{\psi_c\}_{c>0}$ of the
form (\ref{e:diss}). Let $\bbL$ denote the Lebesgue measure on $\bbR$. We can take
$$
\frac{d((\nu\otimes \bbL)\circ \psi_c)}{d(\nu\otimes \bbL)}(y,u) =
\frac{d\nu(y)}{d\nu(y)} \frac{d(u+\ln c)}{du}\equiv 1
$$
as a cocycle for the flow $\{\psi_c\}_{c>0}$. By Lemma 3.1 in Pipiras and Taqqu
\cite{pipiras:taqqu:2002s}, for the dissipative flow $\{\psi_c\}_{c>0}$, the cocycle
$b_c$ is given by
\begin{equation} \label{e:b_c}
b_c(y,u) = \frac{b(\psi_c(y,u))}{b(y,u)}
\end{equation}
with some function $b$ taking values in $\{-1,1\}$. Hence, a {\it 2}-semi-additive
functional $\{j_c\}_{c>0}$ for the flow $\{\psi_c\}_{c>0}$ in (\ref{e:II-semi-add})
satisfies
$$
j_{c_1c_2}(y,u) = c_2^{-(H-1/\alpha)} j_{c_1}(y,u) +  \frac{b(\psi_c(y,u))}{b(y,u)}
j_{c_2}(\psi_{c_1}(y,u))
$$
for all $(y,u)\in Y\times \bbR$ and $c_1,c_2>0$.

To solve this equation for $j_c$, set $\widetilde j_c(y,u) = b(y,u) j_c(y,u)$ so that
$$
\widetilde j_{c_1c_2}(y,u) = c_2^{-(H-1/\alpha)} \widetilde j_{c_1}(y,u) +  \widetilde
j_{c_2}(y,u+\ln c_1).
$$
Substituting $u=0$ into this relation and setting $\ln c_1 = v$ so that $c_1 = e^v$, $c_2
= c$ and $\widetilde j(y,s) = \widetilde j_{e^s}(y,0)$, we obtain that $c_1c_2 = e^{v+\ln
c}$ and
$$
\widetilde j(y,v+\ln c) = c^{-(H-1/\alpha)} \widetilde j(y,v) +  \widetilde j_{c}(y,v).
$$
Hence,
$$
j_c(y,u) = (b(y,u))^{-1} \widetilde j_c(y,u) = \frac{\widetilde j(\psi_c(y,u))}{b(y,u)} -
c^{-(H-1/\alpha)} \frac{\widetilde j(y,u)}{b(y,u)} $$ $$ = b_c(y,u) j(\psi_c(y,u)) -
c^{-(H-1/\alpha)}j(y,u),
$$
by (\ref{e:b_c}), where $j(y,u) = \widetilde j(y,u)/b(y,u)$ is some function.
\end{example}

\begin{example}\label{ex:1-semi-additive-ident-dissip}
In view of (\ref{e:I-semi-add}), a {\it 1}-semi-additive functional $\{g_c\}_{c>0}$ for
the identity flow (\ref{e:identity-flow}) satisfies $g_{c_1c_2}(x)=c_2^{-1} g_{c_1}(x) +
g_{c_2}(x)$ whose solution is
$$
g_c(x) = (c^{-1} - 1)g(x)
$$
for some function $g:X\mapsto \bbR$ (Lemma 3.2 in Pipiras and Taqqu
\cite{pipiras:taqqu:2002s}). The corresponding equation for the dissipative flow
(\ref{e:diss}) is $g_{c_1c_2}(y,u)=c_2^{-1} g_{c_1}(y,u) + g_{c_2}(y,u+\ln c_2)$ where
solution is
$$
g_c(y,u) = g(y,u+\ln c) - c^{-1} g(y,u),
$$
for some function $g:Y\times \bbR \mapsto \bbR$ (Lemma 3.2 in Pipiras and Taqqu
\cite{pipiras:taqqu:2002s}).
\end{example}

\medskip

\noindent{\bf Remark.} Consider Relation (\ref{e:generated}) where the flow
$\{\psi_c\}_{c>0}$ is the identity. Observe that, by using the relations
(\ref{e:jacobian-identity}) and (\ref{e:cocycle-identity}) and Examples
\ref{ex:2-semi-additive} and \ref{ex:1-semi-additive-ident-dissip}, the equation
(\ref{e:generated}) becomes: for any $c>0$,
$$
c^{-(H-1/\alpha)} G(x,cu) = G\Big(x,u+(c^{-1}-1)g(x)\Big) + (1-c^{-(H-1/\alpha)})j(x)
$$
a.e.\ $\mu(dx)du$, when $H\neq 1/\alpha$, and
$$
G(x,cu) = G\Big(x,u+(c^{-1}-1)g(x)\Big) + (\ln c) j(x)
$$
a.e.\ $\mu(dx)du$, when $H= 1/\alpha$. These equations were also obtained in the proof of
Theorem 5.1 in Pipiras and Taqqu \cite{pipiras:taqqu:2002s} and then used to solve for
the function $G$. The arguments in that theorem, leading to the two equations were quite
involved but as we see here, the equations follow easily once {\it 2}-semi-additive
functionals are used.

\begin{example}\label{ex:1-semi-additive-cyclic}
For $v\in\bbR$ and $a>0$, let
\begin{equation}\label{e:fr-int-parts}
[v]_a = \max\{n\in\bbZ : na\leq v\},\quad \{v\}_a = v - a[v]_a.
\end{equation}
By Theorem 2.1 in Pipiras and Taqqu \cite{pipiras:taqqu:2003cy}, any {\it cyclic flow} is
null-isomorphic to the flow
\begin{equation}\label{e:cyclic}
\psi_c(z,v) = (z,\{v+\ln c\}_{q(z)})
\end{equation}
on the space $(X,\mu) = (Z\times [0,q(\cdot)),\sigma(dz)dv)$, where $q(z)>0$ is some
function. Observe that $\{v+\ln c\}_{q(z)}$, as a function of $v$, has the shape of a
seesaw. Relation (\ref{e:I-semi-add}) for the {\it 1}-semi-additive functional
$\{g_c\}_{c>0}$ of the flow (\ref{e:cyclic}) becomes $g_{c_1c_2}(z,v)=c_2^{-1}
g_{c_1}(z,v) + g_{c_2}(z,\{v+\ln c_2\}_{q(z)})$. The solution to this equation, which is
given in Proposition \ref{p:1-semi-add-cyclic} below, is as follows:
\begin{equation}\label{e:1-semi-add-cyclic}
    g_c(z,v) = g(z,\{v+\ln c\}_{q(z)}) - c^{-1}g(z,v),
\end{equation}
for some function $g:Z\times [0,q(\cdot)) \mapsto \bbR$.
\end{example}

\begin{example}\label{ex:2-semi-additive-cyclic}
We now consider {\it 2}-semi-additive functionals for cyclic flows (\ref{e:cyclic}).  By
Lemma 8.2 in Pipiras and Taqqu \cite{pipiras:taqqu:2003in}, the cocycle $b_c(z,v)$ in
(\ref{e:generated}) for the flow $\{\psi_c\}_{c>0}$ can be expressed as
\begin{equation}\label{e:cocycle-cyclic}
b_c(z,v) = b_1(z)^{[v+\ln c]_{q(z)}} \frac{b(\psi_c(z,v))}{b(z,v)}
\end{equation}
for some functions $b_1:Z\mapsto \{-1,1\}$ and $b:Z\times [0,q(\cdot))\mapsto \{-1,1\}$.
The Radon-Nikodym derivatives
$$
\frac{d((\sigma\otimes \bbL)\circ \psi_c)}{d(\sigma\otimes \bbL)}(z,v) \equiv 1
$$
because $d\{v+\ln c_2\}_{q(z)}/dv = 1$ for almost all $v$ since $q(z)$ does not affect
the slope. These Radon-Nikodym derivatives can be taken as a cocycle for the flow
$\{\psi_c\}_{c>0}$. The {\it 2}-semi-additive functional $\{j_c\}_{c>0}$ in
(\ref{e:II-semi-add}) therefore satisfies
$$
j_{c_1c_2}(z,v) = c_2^{-(H-1/\alpha)} j_{c_1}(z,v) + b_1(z)^{[v+\ln c]_{q(z)}}
\frac{b(\psi_c(z,v))}{b(z,v)} j_{c_2}(z,v).
$$
The solution to this equation, which is given in Proposition \ref{p:2-semi-add-cyclic}
below, is as follows:
$$
j_c(z,v) = b_1(z)^{[v+\ln c]_{q(z)}} \frac{b(\psi_c(z,v))}{b(z,v)} j(\psi_c(z,v)) -
c^{-(H-1/\alpha)} j(z,v)
$$
\begin{equation}\label{e:2-semi-add-cyclic}
+ \frac{j_1(z)}{b(z,v)} [v+\ln c]_{q(z)} 1_{\{b_1(z)=1\}}1_{\{H=1/\alpha\}},
\end{equation}
for some functions $j_1:Z\mapsto \bbR$ and $j:Z\times [0,q(\cdot))\mapsto \bbR$.
\end{example}


\section{Semi-additive functionals for cyclic flows}
\label{s:functionals-cyclic}

In this section, we solve the {\it 1}- and {\it 2}-semi-additive functional equations
(\ref{e:I-semi-add}) and (\ref{e:II-semi-add}) for the cyclic flows $\{\psi_c\}_{c>0}$ of
the form (\ref{e:cyclic}). The results are used in Pipiras and Taqqu
\cite{pipiras:taqqu:2003in} to obtain a general form for the kernel $G$ of a mixed moving
average generated by a cyclic flow.

By (\ref{e:fr-int-parts}), one has $\psi_c(z,v) = (z,\{\ln c + v\}_{q(z)}) = (z,\ln
c+v-nq(z))$ when $nq(z)\leq \ln c + v < (n+1)q(z)$ and hence these flows have the special
representation (\ref{e:special-representation}) with $t$ replaced by $\ln c$ and
\begin{equation} \label{e:special-cyclic}
V(z)=z, \quad r_n(z)=n q(z).
\end{equation}
This representation is convenient when applying Corollary \ref{c:solutions-m} to
characterize semi-additive functionals as in the following propositions.

\begin{proposition}\label{p:1-semi-add-cyclic}
Let $\{\psi_c\}_{c>0}$ be a cyclic flow on the space $Z\times [0,q(\cdot))$ given by
(\ref{e:cyclic}), and let $\{g_c\}_{c>0}$ be a {\it 1}-semi-additive functional for the
flow $\{\psi_c\}_{c>0}$ satisfying (\ref{e:I-semi-add}). Then, the solution to the
equation (\ref{e:I-semi-add}) is given by (\ref{e:1-semi-add-cyclic}).
\end{proposition}

\begin{proof}
Example \ref{ex:semi-additive-cocycle2} shows that $J_c(z,v)=c g_c(z,v)$ is a
semi-additive functional in the sense of (\ref{e:semi-additive-cocycle-m}) for the cyclic
flow $\{\psi_c\}_{c>0}$ and the cocycle $\{B_c\}_{c>0}$ defined by $B_c(z,v)=c$. Since
$\psi_c(z,v)$ has a special representation (\ref{e:special-representation}) with $t$
replaced by $\ln c$, and $V$ and $r_n$ defined in (\ref{e:special-cyclic}), Corollary
\ref{c:solutions-m} shows that the semi-additive functional $\{J_c\}_{c>0}$ can be
expressed as the sum of two semi-additive functionals. After substituting $J_c(z,v)=c
g_c(z,v)$ into their expressions, one gets
$$
 g_{c}(z,v)= g_{c}^{(1)} (z,v)+  g_{c}^{(2)}(z,v),
$$
where
\begin{eqnarray*}
  g_c^{(1)}(z,v) &=&  c^{-1} B_c(z,v) g^{(1)}(z,\{v + \ln c\}_{q(z)}) -
  c^{-1} g^{(1)}(z,v), \\
  g_c^{(2)}(z,v) &=& c^{-1} (B_{e^v}(z,0) )^{-1} \sum_{k\in [0,n)}
  B_{e^{r_k(z)}}(z,0) g_1(V^k z),
\end{eqnarray*}
if $r_n(z)\leq \ln c+v < r_{n+1}(z)$, for some measurable functions $g^{(1)}: Z \times
[0,q(\cdot)) \mapsto \bbR$ and $g_1:Z \mapsto \bbR$.

The function $g_c^{(1)}(z,v)$ has the form (\ref{e:1-semi-add-cyclic}) since $c^{-1}
B_c(z,v)=1$. Consider now the function $g_c^{(2)}(z,v)$. Since $r_n(z)=nq(z)$ by
(\ref{e:special-cyclic}), we have $r_n(z)\leq \ln c+v < r_{n+1}(z)$ when $n = [v +\ln
c]_{q(z)}$ using (\ref{e:fr-int-parts}). By using $B_c(z,v)=c$ and
(\ref{e:special-cyclic}), we obtain that
$$
g_c^{(2)}(z,v) = g_1(z) \, e^{-v-\ln c} \sum_{k\in [0,[v+\ln c]_{q(z)})} e^{kq(z)}=
g_0(z)\, e^{-v-\ln c} \Big(e^{[v + \ln c]_{q(z)}q(z)} - 1\Big),
$$
where $g_0(z) = g_1(z)/(e^{q(z)} - 1)$. Applying (\ref{e:fr-int-parts}), we get
$$
g_c^{(2)}(z,v) = g_0(z) (e^{-\{v+\ln c\}_{q(z)}} - c^{-1} e^{-v}),
$$
which has the form (\ref{e:1-semi-add-cyclic}) with $g(z,v) = g_0(z)e^{-v}$. Thus
$g_c^{(1)}+g_c^{(2)}$ has the form (\ref{e:1-semi-add-cyclic}). \ \ $\Box$
\end{proof}

\begin{proposition}\label{p:2-semi-add-cyclic}
Let $\{\psi_c\}_{c>0}$ be a cyclic flow on the space $Z\times [0,q(\cdot))$ given by
(\ref{e:cyclic}). Let also $\{j_c\}_{c>0}$ be a {\it 2}-semi-additive functional for the
flow $\{\psi_c\}_{c>0}$ satisfying (\ref{e:II-semi-add}), and choose a version of the
Radon-Nikodym derivatives
$$
\frac{d((\sigma\otimes \bbL)\circ \psi_c)}{d(\sigma\otimes \bbL)}(z,v) \equiv 1
$$
which is a cocycle for the flow $\{\psi_c\}_{c>0}$, where $\sigma(dz)$ is a measure on
$Z$ and $\bbL$ denotes the Lebesgue measure on $\bbR$. Then, the solution to
(\ref{e:II-semi-add}) is given by (\ref{e:2-semi-add-cyclic}).
\end{proposition}

\begin{proof}
Example \ref{ex:semi-additive-cocycle1} shows that $J_c(z,v) = c^{H-1/\alpha} j_c(z,v)$
is a semi-additive functional related to the cocycle $B_c(z,v) = c^{H-1/\alpha}b_c(z,v)$.
Since $\psi_c(z,v)$ has a special representation (\ref{e:special-representation}) with
$t$ replaced by $\ln c$, and $V$ and $r_n$ given in (\ref{e:special-cyclic}), we can
apply Corollary \ref{c:solutions-m} to express $\{J_c\}_{c>0}$ by the sum of two
functionals. By substituting $J_c(z,v)=c^{H-1/\alpha} j_c(z,v)$ into these expressions,
we get
$$
j_c(z,v) = j_c^{(1)}(z,v) + j_c^{(2)}(z,v),
$$
where
\begin{eqnarray*}
  j_c^{(1)}(z,v) &=&  c^{-(H-1/\alpha)} B_c(z,v) j^{(1)}(z,\{v + \ln c\}_{q(z)}) -
  c^{-(H-1/\alpha)} j^{(1)}(z,v), \\
  j_c^{(2)}(z,v) &=& c^{-(H-1/\alpha)} (B_{e^v}(z,0) )^{-1} \sum_{k\in [0,n)}
  B_{e^{r_k(z)}}(z,0) j_1(V^k z),
\end{eqnarray*}
if $r_n(z)\leq \ln c+v < r_{n+1}(z)$, for some measurable functions $j^{(1)}: Z \times
[0,q(\cdot)) \mapsto \bbR$ and $j_1:Z \mapsto \bbR$.

Since $c^{-(H-1/\alpha)} B_c(z,v) = b_c(z,v)$ and $b_c$ is given by
(\ref{e:cocycle-cyclic}), the function $j_c^{(1)}(z,v)$ has the form of the first two
terms of (\ref{e:2-semi-add-cyclic}). Consider now the function $j_c^{(2)}(z,v)$. Observe
that $r_n(z)\leq \ln c+v < r_{n+1}(z)$ when $n = [v+\ln c]_{q(z)}$, since $r_n(z)=nq(z)$.
Since $B_c(z,v)=c^{H-1/\alpha} b_c(z,v)$, $V^nz=z$ and $r_n(z)=nq(z)$, we obtain that
$$
j_c^{(2)}(z,v) = j_1(z) c^{-(H-1/\alpha)} e^{-(H-1/\alpha) v} (b_{e^v}(z,0) )^{-1}
\sum_{k\in [0,[v+\ln c]_{q(z)})} e^{(H-1/\alpha) kq(z)} b_{e^{kq(z)}}(z,0).
$$
By using the expression (\ref{e:cocycle-cyclic}) of $b_c(z,v)$, we have
$$
b_{e^v}(z,0)= b_1(z)^{[v]_{q(z)}} \frac{b(z,\{v\}_{q(z)})}{b(z,0)}  = \frac{
b(z,v)}{b(z,0)},
$$
since $(z,v)\in Z\times [0,q(\cdot))$ and $b_1\in \{-1,1\}$. We also get that
$b_{e^{kq(z)}}(z,0) = b_1(z)^k$. Hence,
$$
j_c^{(2)}(z,v) = j_0(z) (b(z,v) )^{-1} e^{-(H-1/\alpha) (v+\ln c)} \sum_{k\in [0,[v+\ln
c]_{q(z)})} e^{(H-1/\alpha) kq(z)} b_1(z)^k,
$$
where $j_0(z)=j_1(z)b(z,0)$. If $H \neq 1/\alpha$ or $b_1(z) \neq 1$, then
$e^{(H-1/\alpha) q(z)}b_1(z) \neq 1$. Hence, by using (\ref{e:fr-int-parts}) and
(\ref{e:cocycle-cyclic}),
$$
j_c^{(2)}(z,v) = j(z) (b(z,v))^{-1} e^{-(H-1/\alpha)(v+\ln c)} \Big( e^{(H-1/\alpha))
q(z)[v+\ln c]_{q(z)}} b_1(z)^{[v+\ln c]_{q(z)}} - 1 \Big)
$$
$$
= j(z) \left(\frac{e^{-(H-1/\alpha) \{v+\ln c\}_{q(z)}}}{b(z,v)}b_1(z)^{[v+\ln c]_{q(z)}}
- \frac{e^{-(H-1/\alpha) (v+\ln c)}}{b(z,v)}\right)
$$
$$
= j(z) \left(b_c(z,v) \frac{e^{-(H-1/\alpha) \{v+\ln c\}_{q(z)}}}{b(z,\{v + \ln
c\}_{q(z)})} - c^{-(H-1/\alpha)} \frac{e^{-(H-1/\alpha) v}}{b(z,v)} \right),
$$
where $j(z) = j_0(z)/(e^{(H-1/\alpha) q(z)} b_1(z)-1)$. Thus $j_c^{(2)}(z,v)$ with $H
\neq 1/\alpha$ or $b_1(z) \neq 1$ has the form of the first two terms of
(\ref{e:2-semi-add-cyclic}) (as did $j_c^{(1)}(z,v)$). If $H =1/\alpha$ and $b_1(z)=1$,
then
$$
j_c^{(2)}(z,v) = j_0(z) (b(z,v))^{-1} \sum_{k\in [0,[v + \ln c]_{q(z)})} 1 =
j_0(z)(b(z,v))^{-1} [v + \ln c]_{q(z)},
$$
which has the form of the last term of (\ref{e:2-semi-add-cyclic}). \ \ $\Box$
\end{proof}

\small


\noindent Vladas Pipiras \hfill  Murad S.\ Taqqu

\noindent  Department of Statistics and Operations Research  \hfill Department of
Mathematics and Statistics

\noindent  University of North Carolina at Chapel Hill \hfill Boston University

\noindent  CB\#3260, New West \hfill 111 Cummington St.

\noindent  Chapel Hill, NC 27599, USA \hfill Boston, MA 02215, USA

\noindent {\it pipiras@email.unc.edu} \hfill {\it murad@math.bu.edu}

\end{document}